\documentclass[12pt]{amsart}

\usepackage{fullpage}

\usepackage{amssymb}
\usepackage{amsmath}
\let\newpf\proof \let\proof\relax 
\newenvironment{pf}{\newpf[\proofname]}{\qed\endtrivlist}

\def\ssubset{\Subset}

\def\m{{\operatorname{m}}}

\def\be{\begin{equation}}
\def\ee{\end{equation}}

\def\ba{{\begin{align}}}
\def\ea{{\end{align}}}

\def\0{{\mathbf 0}}

\def\cal{\mathcal}

\newtheorem*{rigidthm}{Rigidity Theorem}
\newtheorem{thm}{Theorem}[section]

\newtheorem{lemma}[thm]{Lemma}

\newtheorem{prop}[thm]{Proposition}

\theoremstyle{remark}
\newtheorem{rem}{Remark}[section]

\numberwithin{equation}{section}

\def \bn {\hfill \\ \smallskip\noindent}

\theoremstyle{definition}

\def\proof{\bn {\bf Proof.} }

\def\note#1
{\marginpar
{\tiny $\leftarrow$
\par
\hfuzz=20pt \hbadness=9000 \hyphenpenalty=-100 \exhyphenpenalty=-100
\pretolerance=-1 \tolerance=9999 \doublehyphendemerits=-100000
\finalhyphendemerits=-100000 \baselineskip=6pt
#1}\hfuzz=1pt}

\newcommand{\di}{\partial}

\newcommand{\ra}{\rightarrow}

\def\ssk{\smallskip}

\def\sm{\setminus}

\newcommand{\inter}{\operatorname{int}}
\renewcommand{\mod}{\operatorname{mod}}
\newcommand{\tl}{\tilde}

\newcommand{\id}{\operatorname{id}}

\newcommand{\Dom}{\operatorname{Dom}}

\newcommand{\de}{{\delta}}

\newcommand{\si}{{\sigma}}
\newcommand{\Om}{{\Omega}}

\renewcommand{\AA}{{\cal A}}

\newcommand{\GG}{{\cal G}}

\newcommand{\MM}{{\cal M}}

\newcommand{\C}{{\mathbb C}}
\newcommand{\D}{{\mathbb D}}

\newcommand{\T}{{\mathbb T}}

\def\B0{{\bold{0}}}

\catcode`\@=12

\def\Empty{}
\newcommand\oplabel[1]{
  \def\OpArg{#1} \ifx \OpArg\Empty {} \else
  	\label{#1}
  \fi}
		
\newcommand{\comm}[1]{}
\newcommand{\comment}[1]{}

\begin{document}

\title[Combinatorial rigidity for unicritical polynomials]
{Combinatorial rigidity\\ for unicritical polynomials}

\author{Artur Avila, Jeremy Kahn, Mikhail Lyubich and Weixiao Shen}

\address{
CNRS UMR 7599, Laboratoire de Probabilit\'es et Mod\`eles al\'eatoires\\
Universit\'e Pierre et Marie Curie--Boite courrier 188\\
75252--Paris Cedex 05, France
}
\email{artur@ccr.jussieu.fr}

 \address{}
\email{kahn@math.sunysb.edu}

\address{Mathematics Department and IMS, Stony Brook University,
    NY 11794, USA}
\email{mlyubich@math.sunysb.edu}

\address{Department of Mathematics, University of Toronto, Ontario Canada
M5S 3G3}
\email{misha@math.toronto.edu}

\address{Mathematics Department,
University of Science and Technology of China, Hefei, 230026,
CHINA}
\email{wxshen@ustc.edu.cn}

\begin{abstract}

We prove that any unicritical polynomial $f_c:z\mapsto z^d+c$
which is  at most finitely renormalizable and has only repelling
periodic points is combinatorially rigid.
It implies that the connectedness locus (the ``Multibrot set'')
is locally connected at the corresponding parameter values.
It generalizes Yoccoz's Theorem for quadratics to the higher degree
case.
\end{abstract}

\setcounter{tocdepth}{1}

\maketitle

\thispagestyle{empty}
\def\IMSmarkvadjust{0 pt}
\def\IMSmarkhadjust{0 pt}
\def\IMSmarkhpadding{0 pt}
\def\IMSpubltext{Published in modified form:}
\def\SBIMSMark#1#2#3{
 \font\SBF=cmss10 at 10 true pt
 \font\SBI=cmssi10 at 10 true pt
 \setbox0=\hbox{\SBF \hbox to \IMSmarkhpadding{\relax}
                Stony Brook IMS Preprint \##1}
 \setbox2=\hbox to \wd0{\hfil \SBI #2}
 \setbox4=\hbox to \wd0{\hfil \SBI #3}
 \setbox6=\hbox to \wd0{\hss
             \vbox{\hsize=\wd0 \parskip=0pt \baselineskip=10 true pt
                   \copy0 \break%
                   \copy2 \break%
                   \copy4 \break}}
 \dimen0=\ht6   \advance\dimen0 by \vsize \advance\dimen0 by 8 true pt
                \advance\dimen0 by -\pagetotal
	        \advance\dimen0 by \IMSmarkvadjust
 \dimen2=\hsize \advance\dimen2 by .25 true in
	        \advance\dimen2 by \IMSmarkhadjust

%
%
  \openin2=publishd.tex
  \ifeof2\setbox0=\hbox to 0pt{}
  \else 
     \setbox0=\hbox to 3.1 true in{
                \vbox to \ht6{\hsize=3 true in \parskip=0pt  \noindent  
                {\SBI \IMSpubltext}\hfil\break
                \input publishd.tex 
                \vfill}}
  \fi
  \closein2
  \ht0=0pt \dp0=0pt
 \ht6=0pt \dp6=0pt
 \setbox8=\vbox to \dimen0{\vfill \hbox to \dimen2{\copy0 \hss \copy6}}
 \ht8=0pt \dp8=0pt \wd8=0pt
 \copy8
 \message{*** Stony Brook IMS Preprint #1, #2. #3 ***}
}

\SBIMSMark{2005/05}{July 2005}{}


\section{Introduction}

Let us consider the one-parameter family of {\it unicritical} polynomials
$$
   f_c: z \mapsto z^d+c,\quad c\in \C,
$$ 
of  degree $d\geq 2$.  
Let $\MM=\MM_d=\{c \in \C,\, \text {the Julia set of } f_c \text { is 
connected}\}$ be the {\it connectedness locus} of this family.  
In the case of quadratic polynomials ($d=2$),
it is known as the Mandelbrot set,
while  in the higher degree case it is sometimes called the
{\it Multibrot set} (see \cite{Sch2}).

{\it Rigidity} is one of the most remarkable phenomena observed in
holomorphic dynamics.
In the unicritical case this phenomenon assumes (conjecturally) a
particularly strong form of   {\it combinatorial rigidity}: 
combinatorially equivalent non-hyperbolic maps are conformally
equivalent. This Rigidity  Conjecture is equivalent to the local
connectivity of the Multibrot sets $\MM_d$. 
In the quadratic case, we are dealing with the
famous MLC conjecture asserting that the Mandelbrot set is locally connected.

About 15 years ago Yoccoz proved that the Mandelbrot set is locally
 connected  at all non-hyperbolic parameter values  
which are at most finitely renormalizable, see \cite {H}.
In fact, this theorem consists of two independent parts
dealing respectively with maps  that 
have neutral  periodic points or not.  
In the presence of neutral points, Yoccoz's method  
 extends readily to the higher degree case.
However, the proof in the absence of neutral points was
linked to the quadratic case in a very significant way.%
\footnote{See also \cite{ALM,K,puzzle,R} for other proofs of this result
  in the quadratic case.}

\begin{rigidthm}

Let $f_c$, $c \in \MM_d$, be an at most finitely renormalizable 
unicritical polynomial with all periodic points repelling.  
Then $f_c$ is combinatorially rigid.

\end{rigidthm}

Our work uses ``complex bounds'' recently proved in \cite {KL2}, which
in turn are  based on new analytic techniques developed in \cite {KL1}.

While combinatorial rigidity is a statement about polynomials with exactly
the same combinatorics in all scales, 
our further analysis (geometric and measure-theoretical)
of the parameter plane \cite{ALS} (with  applications to the
real case) will depend on comparison of polynomials 
whose combinatorics coincide only up to a certain scale. 
For such maps one can consider  {\it pseudo-conjugacies}, that is,
homeomorphisms which are equivariant up to that scale.
In the course of the proof of the Rigidity Theorem, we will show
that these pseudo-conjugacies can be selected uniformly
quasi-conformally, 
generalizing part of the analysis of \cite {puzzle} in the quadratic case.

Let us point out that our argument for existence of pseudo-conjugacies
is considerably simpler than the previous arguments, while
needing much weaker geometric control of the dynamics.  
Also,  though we restrict ourselves to the unicritical case in order
not to overshadow the idea of the method, 
our argument can be  extended to the multicritical case.
 
In conclusion, let us briefly outline the structure of the paper.
In \S \ref{sec 2} we construct a ``favorite nest'' of puzzle pieces and
transfer {\it a priori bounds} of \cite{KL2} to this nest.
The next section, \S \ref{sec 3}, is central in the paper:
here we prove, using the {\it a priori} bounds, that the respective favorite puzzle pieces
of two maps with the same combinatorics stay a bounded Teichm\"uller
distance apart. In the last section, \S \ref{sec 4}, we derive from it,
by means of the ``pullback argument'', our Rigidity  Theorem. 

Note finally that for real polynomials of any degree,
the real version of the Rigidity Theorem 
has been recently proved in \cite{KSS}. 

\medskip
{\bf Basic Notations and terminology.}\\
  $\D_r=\{z \in \C:\, |z|<r\}$, $\D=\D_1$, $\T=\partial \D$;\\
  $\Dom R$ will stand for the domain of a map $R$;\\
  Connected components will be referred to as ``components'';\\
 {\it Pullbacks} of an {\it open} topological disk $V$ under $f$ 
   are components of $f^{-1}(V)$; \\
{\it Pullbacks} of a {\it closed} disk $V$ are the closures of the
 pullbacks of $\inter V$.

\medskip
{\bf Acknowledgment.} This work was supported in part by the NSF and
NSERC.

\section{The complex bounds}\label{sec 2}

In this section we fix a map $f=f_c: z\mapsto z^d+c$.  
The  constants below  may depend implicitly on its degree $d$, but not on $c$.

Let $V$ be an (open or closed) Jordan disk $V\subset \C $.  
We say that $V$ is {\it nice} 
if $f^k(\partial V) \cap \inter V=\emptyset$ for all $k \geq 1$.

Let $R_V: \Dom R_V\ra V$ be the {\it first return map} for $f$ to a nice {\it open} disk $V\ni 0$.
This map has a nice structure:
its  restriction to each component $U$ of $\Dom R_V$ is a proper map onto
$V$. The degree of this restriction is $d$ or $1$ depending on whether $0\in U$ or otherwise.  
In the former case,
 $U$ is called the {\it central component} of $R_V$.

If $V$ is a {\it closed} nice disk with $0\in \inter V$, 
then we can apply the previous discussion to $\inter V$.
Somewhat abusing notations, we will denote $\Dom R_V$ the closure of the
$\Dom R_{\inter V}$ (and we consider $R_V$ only on $\Dom R_V$).
Then the {\it central piece} $W$ of $\Dom R_V$ is defined as the closure 
of the central component $U$ of $\Dom R_{\inter V}$. 
Notice that $W$ is not necessarily a component of $\Dom R_V$. 

The {\it first landing map} $L_V$ to a nice domain $V\ni 0$ has even nicer properties:
it univalently maps each component of $\Dom L_V$ onto $V$
(one of these components is $V$ itself, and $L_V=\id$ on it). 
In case when $V$ is closed, we will apply to the domain of the landing maps
the same conventions as for the return map.

Let us consider two nice disks,  $V$ and $V'$, containing  $0$ in their interior.
We say that $V'$ is a {\it child} of  $V$ 
if there exists $t \geq 1$ such that $f^t:V' \to V$ is a branched 
  covering of degree $d$.
(Note that $V'\subset V$.)
  We can alternatively say that $V$ is the {\it parent} of $V'$ 
(notice that any child has a single parent but not the other way around).

The children of $V$ are naturally ordered by inclusion. 
Notice that the {\it first child} of $V$ coincides with the central piece of $R_V$
(whenever it exists).
We say that $V'$ is a {\it good child} if $f^t(0)\in \inter U$ where $U$ is the
first child of $V$.  In this case, the first child $U'$ of $V'$ is
contained in $(f^t|V')^{-1}(U)$.  In particular
$$\mod(V' \setminus U') \geq \frac {1} {d} \mod(V \setminus U).$$

\medskip
A {\it puzzle} is a graded (by the {\it depth} $k \geq 0$)
collection of nice closed Jordan disks called {\it puzzle pieces},
such that for each $k \geq 0$ the puzzle pieces of depth $k$ have disjoint
interiors, and the puzzle pieces of depth $k+1$ are the pullbacks
of the puzzle pieces of depth $k$ under $f$. 


It may happen that the first child $U$ of $V$ is good: 
then $U$ is called {\it spoiled}.  In this case
$R_V(0)\in U$ and the first return to $V$ is
called {\it central}.

We say that a sequence of nested puzzle pieces
$W^m \subset V^m \subset ... \subset W^0 \subset V^0$ is a
{\it modified principal nest} if
\begin{enumerate}
\item $W^i$ is the first child of $V^i$,
\item $V^{i+1}$ is the oldest unspoiled child of $W^i$. 
 (In other words, $V^{i+1}$ is the first child of $W^i$ if the return
 to $W^i$ is non-central, and is the second child otherwise.)
\end{enumerate}
  See \S 2.2 of \cite{KL2} for a detailed discussion of the
  combinatorics of this nest.

The following is the main technical result of \cite {KL2}:

\begin{thm} \label {KLthm}

There exists $\delta>0$ such that for every $\epsilon>0$ there exists
$n_0>0$ with the following property.
Let $W^m \subset V^m \subset ... \subset W^0 \subset V^0$ be a modified
principal nest for $f$.  If
$\mod(V^0 \setminus W^0)>\epsilon$ and $n_0 \leq n \leq m$
then $\mod(V^n \setminus W^n)>\delta$.

\end{thm}

We will need a slight variation of this result.

Given a nice domain $Q$, let $\m(Q)=\inf \mod(Q \setminus D)$ where the
infimum is taken over all  components $D$
of  $\Dom R_Q$.

\begin{lemma} \label {children}

Let $U$ be the first child of $V$, and let $V'$ be any child of $V$.
Then 
$$ \m(V') \geq \frac {1} {d} \mod(V \setminus U). $$

\end{lemma}

\begin{pf}

Let $k>0$ be such that $f^k(V') = V$.
Given a component $D$ of $\Dom R_{V'}$, let 
$\Om= f^k (D)$.  Notice that $f^j(V') \cap V'=\emptyset$ for $1 \leq j <k$,
so that $\Omega \subset \Dom L_{V'}$.
Let us show that
\be\label{Om} 
  \mod (V\sm \Om)\geq \mod (V\sm U).
\ee
It is obvious if $\Om\subset U$. Otherwise let us consider the
following inclusions:
$$\Om\subset V \cap \Dom L_{V'}\subset V \cap \Dom L_U\subset \Dom R_V\subset V.$$ 
Let $\Om'$ be the component of $\Dom L_U$ containing $\Om$,
and let $\Om''$ be the component of $\Dom R_V$ containing $\Om'$. 
Then
$$\mod(V \setminus \Om) \geq\mod (\Om''\sm \Om') \geq \mod(V \setminus
U),$$
and (\ref{Om}) follows.

 Consequently, 
$$\mod(V' \setminus D) \geq 
\frac {1} {d} \mod(V \setminus \Om) \geq \frac {1} {d} \mod(V \setminus U). $$
\end{pf}

The {\it favorite child} $Q'$ 
of $Q$ is the oldest good unspoiled child of $Q$.
It is constructed as follows. Let $P$ be the first child of $Q$. 
Let $k>0$ be the first moment when  $R^k_Q(0) \in Q \setminus P$, and let $l>0$ be
the first moment when $R^{k+l}_Q(0) \in P$ (so, $k+l$ is the moment of the first
return back to $P$ after first escape from $P$ under iterates of $R_Q$). 
 Then $Q'$ is the pullback of $Q$
under  $R^{k+l}_Q$ that contains $0$. (Compare with the construction
of the domain $\tl A$ in Lemma 2.5   of \cite{KL2}.) 
Note that the first child is never favorite.

\begin{lemma} \label {Y}

Let us consider a nest of four puzzle pieces,
 $P' \subset Q' \subset P \subset Q$, 
such that $P$ is the first child of $Q$, 
$P'$ is the first child of $Q'$, 
and $Q'$ is the favorite child of $Q$. 
 If $V$ is a puzzle piece containing $Q$ and whose first child $U$ is
contained in $Q$ then 
$$\mod(Q' \setminus P') \geq \frac {1} {d^2} \m(V).$$

\end{lemma}

\begin{pf}
  Let the moments $k$ and $l$ have the same meaning as in the above construction
of the favorite child. 
Then  $R^{k+l}_Q|P'$ is a $d$-to-$1$ branched covering onto some
domain 
$D \subset P\subset U$ which is a component of $\Dom L_{Q'}$.  
Hence $$\mod(Q' \setminus P')=\frac {1} {d} \mod (Q \setminus D).$$

Assume $Q' \neq D$.  Then $D$ is contained in a component of the domain
of the first return map to $U$.  Hence by Lemma \ref {children},
$$ 
\mod(Q \setminus D) \geq \mod (U \setminus D) \geq \m(U) \geq \frac
    {1} {d} \m(V),$$
and the conclusion follows. 

Assume now that $Q'=D$, and let $\Om=R_V(Q')$. 
Then $\Om \neq Q$  since $Q'$ is not the first child of $Q$. 
Hence $\Om$ returns to $Q\subset V$ sometime, 
so that it is contained in a component of  $\Dom R_V$.  Thus
$$
\mod(Q \setminus Q') \geq
\mod(U \setminus Q')=\frac
{1} {d} \mod(V \setminus \Om) \geq \frac {1} {d} \m(V),
$$
 and  the result follows.
\end{pf} 

\begin{prop} \label {QP}

There exists $\delta>0$ such that for every $\epsilon>0$ there exists
$n_0>0$ with the following property.
Let $P^m \subset Q^m \subset ... \subset P^0 \subset Q^0$ be a nest of
puzzle pieces such that $P^i$
is the first child of $Q^i$ and $Q^{i+1}$ is the favorite child of $Q^i$. 
If $\mod(Q^0 \setminus P^0)>\epsilon$ and $n_0 \leq n \leq m$
then $\mod(Q^n \setminus P^n)>\delta$.

\end{prop}

\begin{pf}

Let us consider the modified principal nest $W^k \subset V^k
\subset ... \subset W^0 \subset V^0$ which begins with $V^0=Q^0$
and ends at the maximal level  $k$  such that $ V^k\supset Q^m$.  
For any $n=0,1,\dots, m$,
let us define  $ k(n)\in [0, k]$ as the maximal level such that
$V^{k(n)} \supset Q^n$. (In particular, $k=k(m)$ by definition.)   

Let us show that if $n \leq m-2$ then $k(n+2)>k(n)$.
Indeed, since  $Q^{n+1}$ is a child of $Q^n\subset V^{k(n)}$, it is contained in the first
child $W^{k(n)}$ of $V^{k(n)}$.
Since $Q^{n+2}$ is not younger than the second child of
$Q^{n+1}\subset W^{k(n)}$,
it is contained in the second child of $W^{k(n)}$, 
and the latter is contained in $V^{k(n)+1}$. Hence $k(n+2)\geq k(n)+1$.

\smallskip
By Theorem \ref {KLthm} and Lemma \ref {children},
it is enough to show that for every natural $n \in [2, m-2]$ 
we have:
\begin{enumerate}
\item [($*$)] Either $\mod(Q^n \setminus P^n) \geq C^{-1} \m(V^{k(n)})$
or $\mod(Q^{n+1} \setminus P^{n+1}) \geq C^{-1} \m(V^{k(n)})$
\end{enumerate}
for some constant $C>0$ which depends only on $d$. 


\smallskip
If $W^{k(n)} \subset Q^n$, Lemma \ref {Y} yields the latter estimate 
with $C=d^2$. 
So, assume $W^{k(n)}\supset Q^n$.

Let $Z^0=W^{k(n)}$, and let $Z^{i+1}$ be the first child of
$Z^i$.\footnote {So  $Z^0\supset Z^1\supset \dots$ is the {\it principal
nest} that begins with $Z^0$, see \cite {puzzle}.}
If $Z^1 \subset Q^n \subset Z^0$, Lemmas \ref {Y} and \ref {children}
imply that 
$$\mod(Q^{n+1} \setminus P^{n+1}) \geq
\frac {1} {d^2} \m(Z_0) \geq
\frac {1} {d^3} \m(V^{k(n)}), $$
and we are done. 

So, assume  $Q^n \ssubset Z^1$. Let us consider the first return map
$R=R_{Z^0}$, and
let us find the level $j>0$  such that%
\footnote{Thus,  $Z^0\supset Z^1 \supset \dots \supset Z^j$ is a central
  cascade of puzzle pieces, compare \S 2.2 of \cite{KL2}.} 
 $R(0) \in Z^{j-1} \setminus Z^j$. 
 If $j=1$ then $Z^1=V^{k(n)+1}\subset Q^n$,
contradicting the assumption.  So $j>1$.

Let $D \subset Z^0 \setminus Z^1$ be the component of $\Dom R$
containing $R^j (0)$.
Then $V^{k(n)+1}$ is the component of $(R^j)^{-1}(D)$ containing
$0$.
Let $D' \subset D$ be the component of $\Dom L_{Z^j}$
 containing $R^j (0)$.
Then $Z^{j+1}$ is the component of $(R^j)^{-1}(D')$ containing
$0$.  So we have $Z^j \supset V^{k(n)+1} \supset Z^{j+1}$ and
$$ 
\mod(V^{k(n)+1}
\setminus Z^{j+1})=\frac {1} {d} \mod(D \setminus D') \geq \frac {1} {d}
\mod(Z^0 \setminus Z^1).
$$

If $Q^n \subset Z^j$ then $P^n\subset Z^{j+1}$,
and we obtain the following nest:
$$Z^j \supset Q^n \supset V^{k(n)+1} \supset Z^{j+1} \supset P^n.$$
  It follows that
$$
\mod(Q^n \setminus P^n) \geq \mod(V^{k(n)+1} \setminus Z^{j+1}) \geq
\frac {1} {d} \mod (Z^0 \setminus Z^1) \geq \frac {1} {d^2} \m(V^{k(n)})
$$
and we are done.

If $Z^j \subset Q^n \ssubset Z^1$ then
$Q^n$ is the first child of $R(Q^n)$.
Since every child has a single parent,  $R(Q^n)=Q^{n-1}$.  
This is a contradiction since by definition,
$Q^n$ is the favorite (and hence not the first) child of $Q^{n-1}$.
\end{pf}

\section{Teichm\"uller distance between puzzle pieces}\label{sec 3}

A {\it good nest} is a sequence $Q_m \subset ... \subset Q_0$
such that $Q_i$ is a good child of $Q_{i-1}$, $1 \leq i \leq m$.

\begin{thm} \label {pseudoabs}

Let $c,\tilde c \in \C$, and let $f=f_c$, $\tilde f=f_{\tilde c}$.
Let $Q^m \subset ... \subset Q^0$, $\tilde Q^m \subset ... \subset \tilde
Q^0$, be good nests for $f$, $\tilde f$, such
that there exists a homeomorphism $h:\C \to \C$ with $h(Q^i)=\tilde
Q^i$, 
$0\leq i \leq m$, and $h \circ f(x)=\tilde f \circ h(x)$, $x \notin Q^m$.%
\footnote{We refer to this property as {\it combinatorial
    equivalence} of the nests. }
Let $P^i$ and $\tilde P^i$, $0 \leq i \leq m-1$, be the first kids of $Q^i$
and $\tilde Q^i$ respectively.
Assume that
\begin{enumerate}
\item $\mod (Q^i\sm P^i)>\de$ and $\mod (\tilde Q^i\sm \tilde P^i) >\delta$, $0 \leq i \leq m-1$;
\item $h|\partial Q^0$ extends to a $K$-qc map $(Q^0,0) \to (\tilde Q^0,0)$.
\end{enumerate}
Then $h|\partial Q^m$ extends to a $K'$-qc map $(Q^m,0) \to (\tilde Q^m,0)$
where $K'=K'(\delta,K)$.

\end{thm}

The basic step of  the proof of Theorem \ref {pseudoabs} is the following
lemma on covering maps of the disk.

\begin{lemma} \label {covering}

For every $0<\rho <r<1$ there exists $K_0=K_0(\rho,r)$ with the following
property.
Let $g,\tilde g:(\D,0) \to (\D,0)$ be
holomorphic proper maps of degree $d$.
Let $h,h':\T \to \T$, be such that
$\tilde g \circ h'=h \circ g$. Assume that
\begin{enumerate}
\item The critical values of $g$, $\tilde g$ are contained in
$\overline \D_\rho$;
\item $h$ admits a $K$-qc extension $H:\D \to \D$ which is the identity on $\D_{r}$.
\end{enumerate}
Then $h'$ admits a $K'$-qc extension $H':\D \to \D$ which is the
identity on $\D_r$, where $K'=\max \{K,K_0\}$.

\end{lemma}

\begin{pf}

Let  $\GG_\rho$ be the family  of proper holomorphic  maps  
$G :(\D,0) \to (\D,0)$ of degree $d$  whose critical values are contained in
$\overline \D_\rho$, endowed with  the strong
topology of their extensions to rational maps of  degree $d$.
This family is compact. One can see it, e.g.,  by checking normality of this family on the
whole Riemann sphere.
Normality is  obvious on  $\D$ and $\C\setminus \overline\D$.
To see normality near the unit circle $\T$,
notice that the full preimages  $G^{-1}(\D_r\cup (\overline{\C}\sm \D_{1/r}))$, $G \in \GG_\rho$,
contain $0$ and $\infty$, and omit a definite symmetric annulus around $\T$
(of modulus $d^{-1}2\log r$).

For any $G\in \GG_\rho$, the domain $U_G=G^{-1}(\D_r)$ is a Jordan disk with analytic anti-clockwise oriented boundary.
By the  Schwarz Lemma, $\D_r \ssubset U_G$.

Let $G \in \GG_\rho$, and let
$\phi:\partial \D_{r^{1/d}} \to \di U_G$ be an orientation preserving
homeomorphism such that $G \circ \phi(z)=z^d$.  Then $\phi$ is an
analytic diffeomorphism, and there exists a
$L$-qc map $H_\phi:\D \to \D$ such that
$H_\phi|\D_r=\id$  and $H_\phi|\partial \D_{r^{1/d}}=\phi$.
Moreover $L=L(r,\rho)$ by compactness of $\GG_\rho$.

Furthermore,  the given map $H: \D\sm \overline \D_r\ra \D\sm
 \overline \D_r$ lifts to a
 $K$-qc map 
$$
  \hat H:\D \setminus \overline U_g \to \D \setminus \overline
  U_{\tilde g}
$$ 
such that $\hat H|\partial \D= h'$ and $\tilde g
\circ \hat H=H \circ g$ on $\D \setminus \overline U_g$.
By the previous discussion, $\hat H$ extends to a qc map
$H':\D \to \D$ such that $H'|U_g$ is the composition of two $L$-qc maps
which are the identity on $\D_r$.  The result follows with $K_0=L^2$.
\end{pf}

\noindent {\it {Proof of Theorem \ref {pseudoabs}.}}
Let us consider moments $t_i$, $1 \leq i \leq m$,  such that
$f^{t_i}(Q^i)=Q^{i-1}$.
By combinatorial equivalence of our nests, $f^{t_i}(\tl Q^i)=\tl Q^{i-1}$.   
Then  $f^{t_i}: Q^i \to Q^{i-1}$ are proper  holomorphic maps of degree
$d$, and similarly for the second nest.

Let $v_i=f^{t_{i+1}+...+t_m}(0)$,  $0 \leq i \leq m$.  
Since $f^{t_m}(0) \in P^{m-1}$ and $f^{t_i}(P^i) \subset P^{i-1}$,
we have:  $v_i \in P^i$, $0 \leq i \leq m-1$.

Let  $\psi_i:(Q^i,v_i) \to (\D,0)$ be the  uniformizations 
of the domains under consideration by the unit disk,
and  let $g_i=\psi_{i-1} \circ f^{t_i} \circ \psi_i^{-1}$. 
The maps $g_i: (\D,0) \ra (\D,0)$
are unicritical  proper holomorphic  maps of degree $d$.
Let  $u_i=\psi_i(0)$  stand for the critical points of these maps. 

The corresponding objects for the second nest will be marked with
tilde. 

 Let us also consider homeomorphisms  $h_i: \T \to \T$  given by
$h_i=\tilde \psi_i \circ h \circ \psi_i^{-1}$.
They are equivariant with respect to the $g$-actions, i.e.,
 $h_{i-1} \circ g_i=\tilde g_i \circ h_i$.

Let $\psi_i(P^i)=\Om^i$.
Since $\mod(\D \setminus \Om^i) \geq \delta$ and
$\Om^i \ni \psi_i(v_i)=0$,  $0 \leq i \leq m$,
these domains are contained in some
disk $\D_\rho$ with $\rho=\rho(\de)<1$. 
Since $f^{t_i}(0)\in P^{i-1}$,
we conclude that $g_i(u_i)\in \Om^{i-1}\subset \D_\rho$,
 $0 \leq i \leq m$.
The same assertions hold for the second nest. 
So, all the maps $g_i$ and $\tl g_i$  satisfy the assumptions of Lemma \ref{covering}.  

By Assumption  (2) of the theorem we are proving, 
$h_0$ extends to a $K$-qc map $(\D,u_0) \to (\D,\tilde u_0)$. 
 Fix some $r\in (\rho, 1)$.
Since $u_0,\tilde u_0 \in \Om^0\subset  \D_\rho$,
we conclude that  $h_0$ extends to an $L$-qc map $\D \to \D$ which is the
identity on $\D_r$, where
$L=L(K,\delta)$. 

 Let $K_0=K_0(\rho, r)$ be as in Lemma \ref {covering},
and let $K'=\max\{L,K_0\}$.
Consecutive applications of Lemma~\ref {covering} 
show that for $ i=1,\dots, m$, 
the maps  $h_i$ admit a $K'$-qc extension $H_i:\D \to \D$ which are the identity
on $\D_r$.
The desired extension of $h|\partial Q^m$ is now obtained by taking
$\tilde \psi_m^{-1} \circ H_m \circ \psi_m$.
\qed

\section{Pullback Argument}\label{sec 4}

In this section we will derive the Rigidity Theorem from the bound on
the Teichm\"uller distance between the central puzzle pieces by means
of the ``Pullback Argument'' in  the Yoccoz puzzle framework.
This method is standard in holomorphic dynamics.   



\subsection{Combinatorics of a map}\label{comb}

Let $\AA$ stand for the set of parameters $c$ for which the map $f_c:
z\mapsto z^d+c$ has an attracting fixed point. In the quadratic case,
it is a domain bounded by the main cardioid of the Mandelbrot set.
In the higher degree case, $\AA$ is a domain bounded by a simple
closed curve with $d-1$ cusps. 

For the construction of the Yoccoz puzzle for a map $f_c$ with
$c \in \MM \setminus \overline \AA$,
the reader can consult \cite {KL2}, \S 2.3.
Keeping in mind future applications,  here we will extend the construction
(up to a certain depth)
 to some parameters outside $\MM$.

The set $\MM \setminus \overline \AA$ is disconnected.  Each connected
component of $\MM \setminus \overline \AA$ is called a {\it limb}.  The
closure of a limb intersects $\overline \AA$ at a single point called the
{\it root} of the limb.  There are two external rays landing at the root. 
Their union divides $\C$ into two (open) connected components: the one
containing the limb is called a {\it parabolic wake} (see \cite{DH,M2,Sch1}).

For $c$ inside a limb,  the map $f_c$ has a unique 
{\it dividing} repelling fixed point $\alpha$:
the rays landing at it, together with $\alpha$ itself,
disconnect the plane into $q \geq 2$ domains.
This repelling fixed point, and the $q$ external rays landing at it,
have an analytic continuation through the whole parabolic wake.

Let us truncate the parabolic wake by an equipotential of height $h$.
For $c$ in the
truncated parabolic wake, the Yoccoz puzzle pieces of depth $0$ are obtained
by taking the closure of the
connected components of $\C \setminus \overline {\{\text {external rays
landing at } \alpha\}}$ truncated by the equipotential of height $h$.

We say that $f$ {\it has well defined combinatorics up to depth} $n$ if
$$
       f^k(0) \in \bigcup_j \inter Y^0_j,\quad 0 \leq k \leq n.
$$
  In this case we define
Yoccoz puzzle pieces of depth $n$ as the pullbacks
of  the Yoccoz puzzle pieces of depth $0$ under $f^n$.
The puzzle pieces of depth $n$ will be denoted by $Y^n_j$, where the
label $j$ stands for the angles of the external rays that bound $Y^n_j$.
The puzzle piece of depth $n$ whose interior contains $0$ is called
the {\it critical puzzle piece} of depth $n$ and it is also denoted $Y^n$.
The {\it combinatorics of $f$ up to depth} $n$
(provided it is well defined)
 is the set of labels of puzzle pieces of depth $n$.
Note that the combinatorics up to depth $n+t$ determines the puzzle piece $Y^n_j$
containing the critical value $f^t(0)$.

If $f$ does not have well defined combinatorics of all depths, then either
the Julia set of $f$ is disconnected or the critical point is eventually
mapped to the repelling fixed point $\alpha$.  
Otherwise there are critical puzzle pieces of all depth.
In this case, we say that $f$ is {\it combinatorially recurrent}
if the critical point returns to all critical puzzle pieces.
Combinatorially recurrent maps can be either renormalizable or
non-renormalizable, see \cite {KL2}, \S 2.3.

Two non-renormalizable maps are called {\it combinatorially equivalent}
if they have the same combinatorics up to an arbitrary depth.
(See \S \ref{final} for a definition of combinatorial equivalence in the
renormalizable case.)

The following result treats the main special
case of the Rigidity Theorem.


\begin{thm} \label {thmrigi}

Let $f: z\mapsto z^d+c$ be a non-renormalizable combinatorially
 recurrent map.
 If $\tilde f: z\mapsto z^d+\tl c$  is combinatorially equivalent to $f$,
then $f$ and $\tilde f$ are quasiconformally conjugate.

\end{thm}

In the next section we will deduce Theorem \ref {thmrigi}
from a more general statement regarding pseudo-conjugacies.

\subsection{Pseudo-conjugacies and rigidity}

In this section  $f$ will stand for a map  
satisfying  assumptions of  Theorem~\ref{thmrigi}. 
For such a map, 
the construction of the favorite child preceding Lemma \ref{Y}  and 
the discussion of the modified principal nest (see \cite{KL2}, \S 2.2-2.3)
yield:
\begin{enumerate}
\item Every critical puzzle piece $Y^s$ has a favorite child. 

\item Let  $l>0$ be the  minimal moment for which $f^{lq}(0) \notin
 Y^1$.
Then the first child of $Y^{lq}$ is contained in $\inter Y^{lq}$.
\end{enumerate}
This allows us to construct an infinite nest
 $Q^0 \supset P^0 \supset Q^1 \supset P^1 \supset ...$
as follows.  
Take $Q^0= Y^{lq}$,
let $Q^{i+1}$ be the favorite child of $Q^i$, and let $P^i$ be the
first child of $Q^i$.  

If $f$ and $\tilde f$ have the same combinatorics up to depth $n$,
a {\it weak  pseudo-conjugacy} (up to depth $n$) between $f$ and $\tilde f$ is an
orientation preserving homeomorphism $H:(\C,0) \to (\C,0)$ such that
$H(Y^0_j)=\tilde Y^0_j$ and $H \circ f=\tilde f \circ H$
outside the interior of the puzzle pieces of depth $n$.
If the last equation is satisfied everywhere outside the central
puzzle piece $Y^n$, then $H$ is called a {\it pseudo-conjugacy}
(up to depth $n$). 

A (weak) pseudo-conjugacy is said to {\it match the B\"ottcher marking} 
if near $\infty$  it becomes the identity in  the B\"ottcher coordinates 
for $f$ and $\tl f$. 
(Then  by equivariance it is the identity in  the B\"ottcher
coordinates outside $\cup_j Y^n_j$ and $\cup_j Y^n_j$.)  
{\it In what follows all (weak) pseudo-conjugacies are assumed to match the 
B\"ottcher marking.}

The following lemma provides us with a weak pseudo-conjugacy 
(between $f$ and $\tl f$)  with a weak dilatation control. 

\begin{lemma} \label {K_n}
 
If $f$ and $\tilde f$ have the same combinatorics up to depth $n$ then there
exists a $K_n$-qc weak pseudo-conjugacy between
$f$ and $\tilde f$.
(Here $K_n$ depends on the maps $f$ and $\tl f$.)
\end{lemma}

\begin{pf}

The case  $n=0$ can be dealt with by means of holomorphic motions.  
We will only sketch the construction;
details can be found in  \cite {R} (in the case $d=2$ which at
this point does not differ from the higher degree case). 

The property 
that $f$ and $\tilde f$ have the same combinatorics up to
depth $0$ just means that $c$ and $\tilde c$ belong to the same truncated
parabolic wake.  Inside the truncated parabolic wake, the $q$ external rays
landing at the $\alpha$ fixed point, and the equipotential of height $h$,
move holomorphically in $\C \setminus \{0\}$.
Namely, there exists a family of injective maps $\phi_b$,
parametrized by a parameter $b$ in the truncated parabolic wake,
which map the rays and equipotential in question for $c$ 
to the corresponding curves for $b$
(matching the B\"ottcher marking), 
and such that $b \mapsto \phi_b(z)$ is holomorphic, $\phi_c=\id$.

Outside the equipotential of height $h$,
this holomorphic motion extends to a  motion holomorphic in both
variables $(b,z)$ and  tangent to the identity at $\infty$
(it comes from the B\"otcher coordinate near $\infty$).  
By \cite {BR}, the map $\phi_b$ extends to a $K_0$-qc map 
$(\C, 0)\ra (\C,0)$,  where $K_0$  depends only on the
hyperbolic distance between $c$ and $b$ inside the truncated
parabolic wake.  This is the desired qc weak pseudo-conjugacy  $H_0$ for $n=0$.

We will now treat the general case by induction.  Assuming that it holds for some
$n-1 \geq 0$, let us modify the qc weak pseudo-conjugacy $H_{n-1}$ up to depth $n-1$ inside
the puzzle piece of depth $n-1$ containing the critical value $c$, so
that it takes $c$ to $\tilde c$.  The resulting map $H'_{n-1}$ is still a
weak pseudo-conjugacy up to depth $n-1$, and can be taken quasiconformal.  
We now define the desired qc weak  pseudo-conjugacy $H_n$ up to depth $n$ as
the lift of  $H_{n-1}'$ (i.e., $\tilde f \circ H_n=H_{n-1}' \circ  f$)
normalized so that $H_n=H_{n-1}$ near infinity.
\end{pf}

The following lemma gives a two-fold refinement of the previous one:
first, it improves equivariance properties of a weak pseudo-conjugacy $H$ 
turning it into a pseudo-conjugacy $H'$;
more importantly, it provides us with a dilatation control of $H'$
in terms of the Teichm\"uller distance between the deepest puzzle pieces. 

\begin{lemma} \label {Kqc}

Let $H$ be a qc weak pseudo-conjugacy  up to depth $n$
between $f$ and $\tilde f$.
Assume that $H|\partial Y^n$ admits a $K$-qc extension $(\inter Y^n,0) \to
(\inter \tilde Y^n,0)$. 
Then there exists a $K$-qc pseudo-conjugacy (up to depth $n$)
$H'$ between $f$ and $\tilde f$.

\end{lemma}

\begin{pf}

We may assume that $H|\inter Y^n$ is $K$-qc.  Let $H^{(0)}=H$.
Let us  construct by induction a sequence of weak pseudo-conjugacies (up
 to depth $n$)   $H^{(j)}$ as follows. Assume $H^{(j-1)}$ has been
 already constructed.
Since the maps 
$$
   f: \overline \C\sm Y^n \ra \overline \C\sm f(Y^n)\quad {\mathrm and} \quad
   \tl f: \overline \C\sm \tl Y^n \ra \overline \C\sm \tl f(\tl Y^n)
$$
are unbranched coverings of the same degree, the homeomorphism
$$
 H^{(j-1)}:  \overline \C\sm f(Y^n) \ra  \overline \C\sm \tl f (\tl Y^n)
$$
lifts to a homeomorphism
$H^{(j)}:  \overline \C\sm Y^n \ra  \overline \C\sm \tl Y^n$
satisfying the equation 
$H^{(j-1)} \circ f=\tilde f \circ H^{(j)}$
and matching the B\"ottcher coordinate outside the union of  puzzle pieces of
depth $n$. In particular, it matches the B\"ottcher coordinate 
on $\di Y^n$, so it can be extended to $Y^n$ as $H$. 

We obtain a sequence $\{H^{(j)}\}_{j \geq 0}$ of
qc weak pseudo-conjugacies with non-in\-creasing dilatation. 
Hence it is precompact in the uniform topology.
Moreover,   $H^{(j)}=H^{(j-1)}$ 
outside the union of puzzle pieces of depth $n+j-1$.
Thus,  the sequence  $\{H^{(j)}\}$ converges pointwise  outside the
filled Julia set $K(f)$. 
 Since $K(f)$ has empty interior, we conclude that $H^{(j)}$
converges uniformly on the whole plane to some
qc weak pseudo-conjugacy $H'$ up to depth $n$.

By construction, $H'$ coincides with $H$ on $\inter Y^n$ and also
outside the union of puzzle pieces of depth $n$ (in particular it matches
the B\"ottcher marking near $\infty$).  Moreover, $H' \circ f=\tilde f \circ H'$
outside $Y^n$, so that $H'$ is a qc pseudo-conjugacy.
It follows that the dilatation of $H'$ is bounded by $K$
except possibly on the set
$
   X= \{x \in J(f):\, f^k(x) \notin \inter Y^n, \quad k \geq 0\}
$  
(here $J(f)$ stands for the Julia set of $f$).
This set is uniformly expanding, and hence has zero Lebesgue measure.  The result follows.
\end{pf}

\ssk
{\it Remark.} 
One can  construct the above map $H'$ more directly as follows.
First  define $H'$ on the pieces of  $\Dom L_{Y^n}$ as the univalent pullbacks of $H$
(this map is $K$-qc).
Then define $H'$ on $F(f)\sm \Dom L_{Y^n}$ (where $F(f)=\overline \C
\setminus J(f)$ is the Fatou set of $f$)
to be  the identity in the B\"ottcher
coordinates (this map is conformal). These two maps match on the common boundary of the pieces
since $H$ respects the B\"ottcher marking on $\di Y^n$. 
Since the residual set $X$ is hyperbolic,
one can show that this map admits a $K$-qc extension to the whole plane.
\ssk


Let $q_m$ (respectively, $p_m$)  be the depth of the puzzle piece $Q^m$
(respectively, $P^m$), 
i.e.,  $Q^m= Y^{q_m}$ (respectively, $P^m=Y^{p_m}$).

\begin{thm} \label {thmpseudo}

Assume that $f$ is combinatorially recurrent and non-renormalizable.
If $\tilde f$ has the same combinatorics as $f$
up to depth $q_m+p_{m-1}-q_{m-1}$,  then there exists a $K$-qc
pseudo-conjugacy between $f$ and $\tilde f$ up to depth $q_m$,
where $K= K(f,\tl f)$.
\end{thm}

\begin{pf}
For $k=0,\dots, m$,
let $h_k$ be the weak pseudo-conjugacies, up to depth $q_k$, 
constructed in  Lemmas \ref {K_n} and \ref{Kqc} 
(with the weak dilatation control at this moment).

Consider the
sequence of puzzle pieces $\tilde Q^k=h_m(Q^k)=\tilde Y^{q_k}$
for $\tilde f$. 
Let us show that $\tilde Q_k$ is the favorite child of $\tilde Q_{k-1}$ for $1 \leq k \leq m$.  
Indeed, it is clear that $\tilde Q^k$ is a child of $\tilde Q^{k-1}$,
and that this child is not the first. 
Moreover, the combinatorics up to level $p_{k-1}+q_k-q_{k-1}$
determines the puzzle piece of depth $p_{k-1}$ containing 
the critical value of the map $\tilde f^{q_k-q_{k-1}}:\tilde Q^k \to \tilde Q^{k-1}$.
Hence  $\tl f^{q_k-q_{k-1}}(0)\in \inter \tilde P^{k-1}$, so  $\tilde Q^k$ is a good child of $\tilde Q^{k-1}$.  
To see that $\tilde Q^k$ is a favorite child of $\tilde Q^{k-1}$, we
reverse this reasoning to conclude that
for $l\in (q_{k-1},q_k)$, the piece $\tilde Y^l$ cannot be a
good non-spoiled child of $\tilde Q^{k-1}$,
 for otherwise $Y^l$ would be a good non-spoiled child of $Q^{k-1}$.
%

Since $h_m|\partial Q^0=h_0|\partial Q^0$, 
$h_m|\partial Q^0$  extends
to a $K_{q_0}$-qc map $Q^0\ra Q^0$ with  $K_{q_0}=K_{q_0}(f,\tilde f)$.
Moreover, by Proposition \ref {QP}, {\it a priori bounds}
(1) of Theorem~\ref{pseudoabs}
hold for the respective nests of $f$ an $\tl f$. 
Applying this theorem, we conclude that 
$h_m|\partial Q^m$ extends to a $K'$-qc map $Q^m \to \tilde Q^m$, where
$K'=K'(f,\tilde f)$.  The result now follows from
Lemma \ref {Kqc}.
\end{pf}

\begin{rem}

The proof shows that $K(f,\tilde f)$ only depends on $K_{q_0}(f,\tilde f)$,
and on $\mod(Q^0
\setminus P^0)$, $\mod(\tilde Q^0 \setminus \tilde P^0)$.

\end{rem}

\noindent {\it Proof of Theorem \ref {thmrigi}.}
Let $h_n$ be the pseudo-conjugacy up to depth $q_n$ between $f$ and $\tilde
f$  given by Theorem \ref {thmpseudo}.
Since the $h_n$ have uniformly bounded dilatations,
we can take a limit map $h$.  
Then $h$ is a qc map satisfying
$h \circ f=\tilde f \circ h$ outside the filled Julia set $K(f)$.  
Since $K(f)$ has empty interior, $h \circ f=\tilde f \circ h$ 
holds everywhere by continuity.  The result follows.\qed

\subsection{Final Remarks}\label{final}

  The Rigidity Theorem stated in the Introduction is reduced to Theorem \ref{thmrigi} by standard means:

\ssk $\bullet$
  The non-combinatorially recurrent case is simple, and is treated in the same way as in the
quadratic case (see \cite{M1}).

\ssk $\bullet$
Rigidity follows from the qc equivalence of combinatorially equivalent maps by an
open-closed argument.  This argument can be summarized as follows
(see e.g.,  \S 5 of \cite {ideas}).
Combinatorial classes of maps with only
repelling periodic orbits are closed subsets of the parameter plane,
while qc classes are either singletons or open (in one-parameter families)
by the Ahlfors-Bers Theorem.
Thus, if some combinatorial class coincides with a qc class, it must be a singleton.

\ssk $\bullet$
  The case of at most finitely renormalizable maps is reduced to the case of non-renormalizable maps
by means of straightening. Namely, let us consider  two maps 
$f: z\mapsto z^d+c$ and  $\tl f: z\mapsto z^d+\tl c$,
which are exactly $n$ times renormalizable. Then there is a nest of little Multibrot copies,
$$
   \MM\supset\MM^1\supset\dots\supset \MM^n\ni \{c,\tl c\},
$$   
such that under the straightening $\si: \MM^n\ra\MM$ the parameters $c$ and $\tl c$
become non-renormalizable. We say that $f$ and $\tl f$ are combinatorially equivalent
if  the corresponding non-renormalizable maps $z\mapsto z^d+\si(c)$
and $z\mapsto z^d+\si(\tl c)$ are%
\footnote{Two infinitely renormalizable maps are called {\it combinatorially equivalent}
  if they belong to the same infinite nest of little Multibrot copies.}
(see discussion in  \cite{Sch2}). 
If so then by the non-renormalizable case of the
Rigidity Theorem, $\si(c)=\si(\tl c)$, and we are done.

\end{document}